\documentclass{birkjour}
\usepackage{amsthm}
\usepackage{mathtools}
\usepackage{url}
\usepackage{helvet}
\usepackage{times}
\usepackage{amsmath}
\usepackage{amsfonts}
\usepackage{amssymb}
\usepackage{mathrsfs}
\usepackage{bbm}
\usepackage{floatflt}
\usepackage{graphics}
 \newtheorem{thm}{Theorem}[section]

 \newtheorem{prop}[thm]{Proposition}
 \theoremstyle{definition}
 
 \theoremstyle{remark}

 \numberwithin{equation}{section}
\def\cubic{{M}^{(3)}_n}
\begin{document}
\title{Quantum super Nambu bracket of cubic supermatrices and 3-Lie superalgebra}
\author{Viktor Abramov}
\address{Institute of Mathematics, University of Tartu\\
Liivi 2 -- 602, Tartu 50409, Estonia}
\email{viktor.abramov@ut.ee}


\date{}
%----------additions
%\dedicatory{To my boss}
%%% ----------------------------------------------------------------------

\begin{abstract}
We construct the graded triple Lie commutator of cubic supermatrices, which we call the quantum super Nambu bracket of cubic supermatrices, and prove that it satisfies the graded Filippov-Jacobi identity of $3$-Lie superalgebra. For this purpose we use the basic notions of the calculus of 3-dimensional matrices, define the $\mathbb Z_2$-graded (or super) structure of a cubic matrix relative to one of the directions of a cubic matrix and the super trace of a cubic supermatrix. Making use of the super trace of a cubic supermatrix we introduce the triple product of cubic supermatrices and find the identities for this triple product, where one of them can be regarded as the analog of ternary associativity. We also show that given a Lie algebra one can construct the $n$-ary Lie bracket by means of an $(n-2)$-cochain of given Lie algebra and find the conditions under which this $n$-ary bracket satisfies the Filippov-Jacobi identity, thereby inducing the structure of $n$-Lie algebra. We extend this approach to $n$-Lie superalgebras.
\end{abstract}

%%% ---------------------------------------------------------------------
\maketitle
%#################################################
%#################################################
\section{Introduction}
$n$-Lie algebra, where $n\geq 2$, is a generalization of the notion of Lie algebra. An integer $n$ indicates the number of elements of algebra necessary to form a Lie bracket. Thus, a Lie bracket of $n$-Lie algebra $\frak g$ is an $n$-ary multilinear mapping $\frak g\times\frak g\times\ldots\times\frak g\;(n\; \mbox{times})\to\frak g$, which is skew-symmetric and satisfies the Filippov-Jacobi identity (also called fundamental identity).

\vskip.3cm
\noindent
A notion of $n$-Lie algebra was proposed by V.T. Filippov in \cite{Filippov}. Earlier Y. Nambu in \cite{Nambu} proposed a generalization of Hamiltonian mechanics by means of a ternary (or, more generally, $n$-ary) bracket of functions determined on a phase space. This $n$-ary bracket can be regarded as an analog of the Poisson bracket in Hamiltonian mechanics and later this $n$-ary bracket was called the Nambu bracket. The dynamics of this generalization of Hamiltonian mechanics is based on the Nambu-Hamilton equation of motion that contains $n-1$ Hamilton functions. It is worth to mention that Y. Nambu proposed and developed his generalization of Hamiltonian mechanics based on a notion of triple bracket with the goal to apply this approach to quarks model, where baryons are combinations of three quarks. Later it was shown that the Nambu bracket satisfied the Filippov-Jacobi identity and thus, the space of functions endowed with the $n$-ary Nambu bracket can be regarded as the example of $n$-Lie algebra. Further development this direction has received in the paper \cite{Takhtajan}, where the author introduced a notion of Nambu-Poisson manifold, which can be regarded as an analog of the notion of Poisson manifold in Hamiltonian mechanics. In the same paper the author also treated the problem of quantization of Nambu bracket and constructed the explicit representation of Nambu-Heisenberg commutation relations by means of a primitive cubic root of unity.

\vskip.3cm
\noindent
Another area of theoretical physics, where the 3-Lie algebras are applied, is a field theory. Particularly the authors of the paper \cite{Basu-Harvey} proposes a generalization of the Nahm's equation by means of quantum Nambu 4-bracket and show that their generalization of the Nahm's equation describes M2 branes ending on M5 branes.

\vskip.3cm
\noindent
The notion of 3-Lie algebra can be extended to a $\mathbb Z_2$-graded case by means of graded Lie bracket and the graded Filippov-Jacobi identity. The classification of simple linearly compact $n$-Lie superalgebras is given in \cite{Cantarini-Kac}.

\vskip.3cm
\noindent
The question of quantization of the Nambu bracket was considered in a number of papers and the first step in this direction was taken already by Y. Nambu, but so far this is the outstanding problem. In the paper \cite{Awata-Li-Minic-Yaneya} the authors propose the realization of quantum Nambu bracket given by the formula
\begin{equation}
[X,Y,Z]=\mbox{Tr}\,X\;[Y,Z]+\mbox{Tr}\,Y\;[Z,X]+\mbox{Tr}\,Z\;[X,Y],
\label{formula 1 in introduction}
\end{equation}
where $X,Y,Z$ are either square matrices or cubic matrices, and prove that this quantum Nambu bracket satisfies the Filippov-Jacobi identity. Inspired by (\ref{formula 1 in introduction}) in \cite{Abramov2017} we extend this approach to $\mathbb Z_2$-graded case by assuming that $X,Y,Z$ are supermatrices and taking supertrace of a matrix instead of the trace and $\mathbb Z_2$-graded commutator instead of ordinary commutator, i.e. we define the graded triple commutator of supermatrices by
\begin{equation}
[X,Y,Z]=\mbox{Str}\,X\;[Y,Z]+(-1)^{x\,\overline{yz}}\mbox{Str}\,Y\;[Z,X]+(-1)^{z\,\overline{xy}}\mbox{Str}\,Z\;[X,Y],
\label{formula 2 in introduction}
\end{equation}
where $X,Y,Z$ are square supermatrices, $x,y,z$ are the degrees of supermatrices, $\overline{xy}=x+y$ and $[\;,\;]$ is the graded commutator of two supermatrices. In analogy with (\ref{formula 1 in introduction}) the graded triple commutator (\ref{formula 2 in introduction}) is referred to as the \emph{quantum super Nambu bracket of square matrices}. We prove that this graded triple commutator satisfies the graded Filippov-Jacobi identity. Hence the super vector space of supermatrices equipped with the graded triple commutator (\ref{formula 2 in introduction}) becomes the 3-Lie superalgebra.

\vskip.3cm
\noindent
The aim of this paper is to extend (\ref{formula 2 in introduction}) to 3-dimensional supermatrices. For this purpose we define a $\mathbb Z_2$-graded structure of a cubic matrix and call a matrix endowed with this structure a \emph{cubic supermatrix}. We define the \emph{super trace of a cubic supermatrix} relative to one of three directions of a cubic matrix. Then we consider the binary product of two cubic matrices, which is formulated by means of products of sections of certain orientation of cubic matrices. Making use of this product we construct with the help of super trace of a cubic supermatrix the triple product of cubic supermatrices and find the identities for this triple product. Then we apply this triple product to quantum super Nambu bracket and obtain the graded triple commutator of cubic supermatrices, which we call \emph{quantum super Nambu bracket of cubic supermatrices}. We prove that quantum super Nambu bracket of cubic supermatrices satisfies the graded Filippov-Jacobi identity and thus induces the structure of 3-Lie superalgebra on the super vector space of cubic supermatrices. In our approach we use the notions and structures of the calculus of 3-dimensional matrices developed in \cite{Sokolov}. We also generalize (\ref{formula 1 in introduction}) by means of cochains of Lie algebra. Given a Lie algebra we construct the $n$-ary Lie bracket by means of a $(n-2)$-cochain of given Lie algebra and find the conditions under which this $n$-ary bracket satisfies the Filippov-Jacobi identity, thereby inducing the structure of $n$-Lie algebra. Particular case of this approach applied to matrix Lie algebra of $n$th order matrices $\frak gl(n)$ is the matrix 3-Lie algebra with quantum Nambu bracket, which is defined with the help of the trace of a matrix (\ref{formula 1 in introduction}). Similarly we propose a more general approach than given by (\ref{formula 2 in introduction}) by making use of cochains of Lie superalgebra and find the conditions under which a cochain induces the structure of $n$-Lie superalgebra.
%*******************************
%********************************
%%%%%%%%%%%%%%%%%%%%%%%%%%%%%%%%%%%%%%%%%%%%%%%%%%%%%%%%%%%%%%%%%%%%%%%%%%
%##########################################################################
%#########################################################################
%################# NEW SECTION
\section{Quantum Nambu bracket and matrix 3-Lie algebra}
In \cite{Awata-Li-Minic-Yaneya} the authors show that given a matrix Lie algebra $\frak{gl}(n)$ of $n$th order square matrices one can construct the 3-Lie algebra by means of the trace of a matrix. In this section we generalize the approach proposed in \cite{Awata-Li-Minic-Yaneya} by means of cochains of Chevalley-Eilenberg complex of a Lie algebra.

\vskip.3cm
\noindent
A concept of \emph{$n$-Lie algebra}, where $n\geq 2$ is an integer, was independently introduced by V.T. Filippov \cite{Filippov} and Y. Nambu \cite{Nambu}. We remind that an $n$-Lie algebra is a vector space endowed with an $n$-ary Lie bracket
$$
(x_1,x_2,\ldots,x_n)\in \frak g\times\frak g\times\ldots\times\frak g\;(n\;\mbox{times})\to [x_1,x_2,\ldots,x_n]\in \frak g,
$$
where an $n$-ary Lie bracket is skew-symmetric and satisfies the Filippov-Jacobi or the fundamental identity (FI)
\begin{eqnarray}
[x_1,x_2,\!\!\!&\ldots&\!\!\!,x_{n-1},[y_1,y_2,\ldots,y_n]] =\nonumber\\
    &&\sum_{k=1}^n[y_1,y_2,\ldots,y_{k-1},[x_1,x_2,\ldots,x_{n-1},y_k],y_{k+1},\ldots,y_n].
   \label{n-ary Filippov identity}
\end{eqnarray}
Particularly a \emph{3-Lie algebra} is a vector space $\frak g$ equipped with a \emph{ternary Lie bracket} $[\;,\;,\;]:\frak g\times\frak g\times\frak g\to \frak g$, which satisfies
\begin{eqnarray}
[a,b,c]=[b,c,a]=[c,a,b],\;\;\;[a,b,c]=-[b,a,c]=-[a,c,b]=-[c,b,a],
\label{symmetries of ternary commutator 2}
\end{eqnarray}
and the \emph{ternary Filippov-Jacobi identity}
\begin{eqnarray}
 [a,b[c,d,e]] = [[a,b,c],d,e]+[c,[a,b,d],e]+[c,d,[a,b,e]].
\label{ternary Filippov identity}
\end{eqnarray}
Given the matrix Lie algebra $\frak{gl}(n)$ of $n$th order matrices one can introduce the ternary Lie bracket of $n$th order matrices as follows
\begin{equation}
[A,B,C]=\mbox{Tr}\,A\;[B,C]+\mbox{Tr}\,B\;[C,A]+\mbox{Tr}\,C\;[A,B].
\label{triple commutator of square matrices}
\end{equation}
Evidently this ternary Lie bracket obeys the symmetries (\ref{symmetries of ternary commutator 2}). It is proved in \cite{Awata-Li-Minic-Yaneya} that this ternary Lie bracket satisfies the ternary Filippov-Jacobi identity and the authors propose to call this ternary Lie bracket \emph{quantum Nambu bracket}. Hence the vector space of $n$th order matrices equipped with the quantum Nambu bracket (\ref{triple commutator of square matrices}) becomes the matrix 3-Lie algebra.

\vskip.3cm
\noindent
Now our aim is to extend this approach to a more general case than the trace of a matrix. Let $\frak g$ be a finite dimensional Lie algebra and $\frak g^\ast$ be its dual. Fix an element of the dual space $\omega\in \frak g^\ast$ and by analogy with (\ref{triple commutator of square matrices}) define the triple product as follows
\begin{equation}
[x,y,z]=\omega(x)\,[y,z]+\omega(y)\,[z,x]+\omega(z)\,[x,y], \quad x,y,z\in \frak g.
\label{ternary Lie bracket with omega}
\end{equation}
Obviously this triple product is symmetric with respect to cyclic permutations of $x,y,z$ and skew-symmetric with respect to non-cyclic permutations, i.e. it obeys the symmetries (\ref{symmetries of ternary commutator 2}). Straightforward computation of the left hand side and the right hand side of the Filippov-Jacobi identity
\begin{eqnarray}
 [x,y[u,v,t]] = [[x,y,u],v,t]+[u,[x,y,v],t]+[u,v,[x,y,t]],
\label{ternary Filippov identity 1}
\end{eqnarray}
shows that one can split all the terms into three groups, where two of them vanish because of binary Jacobi identity and skew-symmetry of binary Lie bracket. The third group can be split into subgroups of terms, where each subgroup has the same structure and it is determined by one of the commutators $[x,y],[u,v],[u,t],[v,t]$. For instance, if we collect all the terms containing the commutator $[x,y]$ then we get the expression
\begin{equation}
(\omega(u)\,\omega([v,t])+\omega(v)\,\omega([t,u])+\omega(t)\,\omega([u,v]))\,[x,y].
\label{sum of three terms}
\end{equation}
Hence the triple product (\ref{ternary Lie bracket with omega}) will satisfy the ternary Filippov-Jacobi identity if for any elements $u,v,t\in \frak g$ we require
\begin{equation}
\omega(u)\,\omega([v,t])+\omega(v)\,\omega([t,u])+\omega(t)\,\omega([u,v])=0.
\end{equation}
Now we consider $\omega$ as a $\mathbb C$-valued cochain of degree one of the Chevalley-Eilenberg complex of a Lie algebra $\frak g$. Making use of the coboundary operator $\delta:\wedge^k\frak g^\ast\to\wedge^{k+1}\frak g^\ast$ we obtain the degree two cochain $\delta\omega$, where $\delta\omega(x,y)=\omega([x,y])$. Finally we can form the wedge product of two cochains $\omega\wedge\delta\omega$, which is the cochain of degree three and
$$
\omega\wedge\delta\omega(u,v,t)=\omega(u)\,\omega([v,t])+\omega(v)\,\omega([t,u])+\omega(t)\,\omega([u,v]).
$$
Hence the third group of terms (\ref{sum of three terms}) of the Filippov-Jacobi identity vanishes if $\omega\in \frak g^\ast$ satisfies
\begin{equation}
\omega\wedge\delta\omega=0.
\label{condition for omega}
\end{equation}
Thus if an 1-cochain $\omega$ satisfies the equation (\ref{condition for omega}) then the triple product (\ref{ternary Lie bracket with omega}) is the ternary Lie bracket and by analogy with (\ref{triple commutator of square matrices}) we will call this ternary Lie bracket \emph{the quantum Nambu bracket induced by a 1-cochain}.

\vskip.3cm
\noindent
We can generalize the quantum Nambu bracket (\ref{ternary Lie bracket with omega}) if we consider a $\mathbb C$-valued 1-cochain of degree $n-2$, i.e. $\omega\in\wedge^{(n-2)}\frak g^\ast$, and define the $n$-ary product as follows
\begin{equation}
[x_1,x_2,\ldots,x_n]=\sum_{i<j}\,(-1)^{i+j+1}\;\omega(x_1,x_2,\ldots,\hat x_i,\ldots,\hat x_j,\ldots,x_n)\;[x_i,x_j].
\label{n-ary quantum Nambu bracket}
\end{equation}
\begin{thm}
Let $\frak g$ be a finite dimensional Lie algebra, $\frak g^\ast$ be its dual and $\omega$ be a cochain of degree $n-2$, i.e. $\omega\in\wedge^{n-2}\frak g^\ast$. The vector space of Lie algebra $\frak g$ equipped with the $n$-ary product (\ref{n-ary quantum Nambu bracket})
is the $n$-Lie algebra, i.e. the $n$-ary product (\ref{n-ary quantum Nambu bracket}) satisfies the Filippov-Jacobi identity (\ref{n-ary Filippov identity}) if and only if a $(n-2)$-cochain $\omega$ satisfies
\begin{equation}
\omega\wedge\delta\omega=0.
\label{condition for omega in the theorem}
\end{equation}
Particularly the vector space of a Lie algebra $\frak g$ endowed with the $n$-ary Lie bracket (\ref{n-ary quantum Nambu bracket}) is the $n$-Lie algebra if $\omega$ is an $(n-2)$-cocycle, i.e. $\delta\omega=0$.
\label{theorem 1}
\end{thm}
\noindent
Particularly from this theorem it immediately follows that the quantum Nambu bracket for matrices  (\ref{triple commutator of square matrices}) satisfies the ternary Filippov-Jacobi identity. Indeed in this case $\omega=\mbox{Tr}$ and this is 1-cocycle because for any two matrices $A,B$ we have
$$\delta\mbox{Tr}\,(A,B)=\mbox{Tr}\,([A,B])=0.$$
In \cite{Arnlind-Makhlouf-Silvestrov} it is proposed to call the 3-Lie algebra constructed by means of an analog of trace the 3-Lie algebra induced by a Lie algebra. In our approach we will use the similar terminology and call the $n$-Lie algebra constructed by means of a $(n-2)$-cochain $\omega$ (\ref{n-ary quantum Nambu bracket}) the \emph{$n$-Lie algebra induced by a Lie algebra $\frak g$ and an $(n-2)$-cochain} $\omega.$ The $n$-ary Lie bracket will be referred to as the \emph{$n$-ary quantum Nambu bracket induced by a $(n-2)$-cochain}.
%%%%%%%%%%%%%%%%%%%%%%%%%%%%%%%%%%%%%%%%%%%%%%%%%%%%%%%%%%%%%
%%%%%%%%%%%%%%%%%%%%%%%%%%%%%%%%%%%%%%%%%%%%%%%%%%%%%%%%%%%%%
%############################################################
%############################################################
% NEW SECTION
\section{Quantum super Nambu bracket and 3-Lie superalgebra}
%************************************************************
The aim of this section is to extend the approach developed in the previous section to 3-Lie superalgebras.

\vskip.3cm
\noindent
Let $\frak g=\frak g_0\oplus\frak g_1$ be a super vector space. We will denote the degree of a homogeneous element $x$ by $\bar x$. Thus $\bar x\in {\mathbb Z}_2$. A super vector space $\frak g$ is an \emph{$n$-Lie superalgebra} if it is endowed with an $n$-ary graded Lie bracket
$$
(x_1,x_2,\ldots,x_n)\in \frak g\times\frak g\times\ldots\times\frak g\;(n\; \mbox{factors}\,)\to [x_1,x_2,\ldots,x_n]\in\frak g,
$$
which is graded skew-symmetric and satisfies the graded Filippov-Jacobi identity or the graded fundamental identity (GFI)
\begin{eqnarray}
&&[x_1,x_2,\ldots,x_{n-1},[y_1,y_2,\ldots,y_n]]=\nonumber\\ &&\;\;\;\;\sum_{k=1}^n (-1)^{\alpha_k}
             [y_1,y_2,\ldots,y_{k-1},[x_1,x_2,\ldots,x_{n-1},y_k],y_{k+1},\ldots,y_n],
\end{eqnarray}
where
$$
\alpha_k=\sum_{i=1}^{n-1} \bar x_i\;\sum_{j=1}^{k-1} \bar y_j.
$$
Particularly a super vector space $\frak g =\frak g_0\oplus\frak g_1$ is said to be a 3-Lie superalgebra if it is endowed with a graded triple Lie bracket $[\;.\;,\;.\;,\;.\;]:\frak g\times\frak g\times \frak g\to \frak g$, i.e.
$$
[\mathfrak g_\alpha,\frak g_\beta,\frak g_\gamma]\subseteq \frak g_{\alpha+\beta+\gamma},\;
     \alpha,\beta,\gamma\in {\mathbb Z}_2,
$$
which is graded skew-symmetric
\begin{equation}
[x,y,z] = -(-1)^{\bar x\bar y}[y,x,z],\label{skew-symmetry 1}
\end{equation}
\begin{equation}
[x,y,z] = -(-1)^{\bar y\bar z}[x,z,y],\label{skew-symmetry 2}
\end{equation}
\begin{equation}
[x,y,z] = -(-1)^{\bar x\bar y+\bar y\bar z+\bar x\bar z}[z,y,x],\label{skew-symmetry 3}
\end{equation}
and satisfies the ternary graded Filippov-Jacobi identity
\begin{eqnarray}
&&[x,y,[z,v,w]]] = [[x,y,z],v,w]+(-1)^{\overline{xy}\,\bar z}[z,[x,y,v],w]\nonumber\\
       &&\qquad\qquad\qquad\qquad\qquad\qquad\;\;\;\qquad +(-1)^{\overline{xy}\,\overline{zv}}[z,v,[x,y,w]].
\label{graded Filippov identity}
\end{eqnarray}

\vskip.3cm
\noindent
Assume $\frak g=\frak g_0\oplus\frak g_1$ is a Lie superalgebra. Now we extend our formula for the quantum Nambu bracket induced by 1-cochain (\ref{ternary Lie bracket with omega}) to the case of a Lie superalgebra. We take $\omega\in\frak g^\ast$ and define the triple product of three elements of a Lie superalgebra by
\begin{equation}
[x,y,z]=\omega(x)\,[y,z]+(-1)^{\bar x\,\overline{yz}}\,\omega(y)\;[z,x]+(-1)^{\bar z\,\overline{xy}}\,\omega(z)\;[x,y],
\label{graded triple product with omega}
\end{equation}
where $[\;,\;]$ is the graded Lie bracket of a Lie superalgebra $\frak g$. It is easy to verify that the triple product is graded skew-symmetric, i.e. it satisfies the symmetries (\ref{skew-symmetry 1}) -- (\ref{skew-symmetry 3}).

\vskip.3cm
\noindent
Now we consider the ternary graded Filippov-Jacobi identity (\ref{graded Filippov identity}) in the case, where a graded skew-symmetric product is defined by the formula (\ref{graded triple product with omega}). By analogy with the case of the Filippov-Jacobi identity for a 3-Lie algebra considered in the previous section we might expect that there will be groups of terms which vanish by virtue of the binary graded Jacobi identity, pairs of terms which cancel by virtue of graded skew-symmetry of the triple product and finally groups of terms which will vanish if we impose a condition on a cochain $\omega$. The simplest case is when the pairs of terms with opposite signs (which arise by virtue of graded skew-symmetry) cancel, no additional conditions are required here. But in the case of terms which can be collected to the graded (binary) Jacobi identity, we have a situation different from the one considered in the previous section. For instance, at the left hand side of the graded Filippov-Jacobi identity we have the term
$$
\omega(x)\,\omega(z)\;[y,[v,w]].
$$
The terms on the right hand side of the graded Filippov-Jacobi identity which could complete the graded (binary) Jacobi identity for this term are
\begin{equation}
\omega(x)\,\omega(z)\;\big((-1)^{\bar z\,\overline{xy}}\;[[y,v],w]+(-1)^{\overline{xy}\,\overline{zv}}\;[v,[y,w]]\big).
\label{expression 1}
\end{equation}
Consequently we get the binary graded Jacobi identity (and then all three terms will vanish) if inside the parenthesis we have the expression
\begin{equation}
[[y,v],w]+(-1)^{\overline{yv}}\;[v,[y,w]].
\label{expression 2}
\end{equation}
It is easy to see that the expression in parenthesis in (\ref{expression 1}) gives the expression (\ref{expression 2}) only in one case, when $x,z$ are even elements, i.e. $\bar x=\bar z=\bar 0.$ Indeed then $\bar z\,\overline{xy}=\bar 0$ and $\overline{xy}\,\overline{zv}=\bar y\,\bar v$. In all other cases we have the improper coefficients in the expression (\ref{expression 1}). Consequently the product $\omega(x)\,\omega(z)$ must be zero as soon as at least one of elements $x,y$ is odd. This can be achieved if for any odd element $x$ of a Lie superalgebra we require $\omega(x)=0$. Hence the first condition for a cochain $\omega$ can be written in the form
\begin{equation}
\omega|_{\frak g_1}\equiv 0.
\label{condition 1 for omega}
\end{equation}
The remained terms in the graded Filippov-Jacobi identity are the following: at the left hand side there is the expression $\omega([z,v,w])\;[x,y]$ and the similar expressions at the right hand side are $\omega([x,y,z])\;[v,w], \omega([x,y,v])\;[z,w], \omega([x,y,w])\;[z,v]$ (all expressions are given up to a sign). The only way to get rid of these terms is to impose the condition
\begin{equation}
\omega([x,y,z])=0,
\label{condition 2 for omega}
\end{equation}
where $x,y,z$ are either all even elements or two of them are odd and one is even. Substituting the expression (\ref{graded triple product with omega}) into the condition (\ref{condition 2 for omega}) we can write this condition in the form
\begin{equation}
\omega(x)\,\omega([y,z])+(-1)^{\bar x\,\overline{yz}}\omega(y)\,\omega([z,x])+(-1)^{\bar z\,\overline{xy}}\omega(z)\,\omega([x,y])=0.
\label{condition 2prime for omega}
\end{equation}
This condition simplifies in the case of two odd degree elements and takes on the form
$$
\omega(z)\,\omega([x,y])=0,\quad z\in\frak g_0,\;x,y\in\frak g_1.
$$
Now $\delta\omega(x,y)=\omega([x,y])$ is the 2-cochain of a Lie superalgebra $\frak g$ and it is the element of the superexterior algebra $\wedge\frak g^\ast$. From the previous considerations it follows that $\delta\omega(x,y)=0$ if either $x\in \frak g_1,y\in\frak g_0$ or $x\in \frak g_0, y\in\frak g_1$. The condition (\ref{condition 2prime for omega}) can be written in the form
\begin{equation}
\omega\wedge\delta\omega=0,
\end{equation}
where we use the wedge product of superexterior algebra $\wedge \frak g^\ast$.
\begin{thm}
Let $\frak g=\frak g_0\oplus\frak g_1$ be a Lie superalgebra, $\frak g^\ast$ be its dual and $\omega\in \frak g^\ast$ be a 1-cochain. Define the graded triple product of elements of a Lie superalgebra $\frak g$ by
\begin{equation}
[x,y,z]=\omega(x)\,[y,z]+(-1)^{\bar x\,\overline{yz}}\,\omega(y)\;[z,x]+(-1)^{\bar z\,\overline{xy}}\,\omega(z)\;[x,y].
\label{graded triple Lie bracket with omega}
\end{equation}
This graded triple product satisfies the graded Filippov-Jacobi identity (graded fundamental identity) if and only if
\begin{enumerate}
\item[i)]
$\omega|_{\frak g_1}\equiv 0$,
\item[ii)]
$\omega\wedge\delta\omega=0$, where $\delta\omega(x,y)=\omega([x,y])$ is 2-cochain and $\wedge$ is the wdege product of superexterior algebra of $\frak g^\ast$.
\end{enumerate}
Thus the super vector space $\frak g$ endowed with the graded triple Lie bracket (\ref{graded triple Lie bracket with omega}) is the 3-Lie superalgebra.
\label{theorem about 3-Lie superalgebra}
\end{thm}
\noindent
In analogy with the terminology of the previous section the 3-Lie superalgebra constructed by means of the graded triple Lie bracket (\ref{graded triple Lie bracket with omega}) will be referred to as \emph{the 3-Lie superalgebra induced by a Lie superalgebra $\frak g$ and a cochain $\omega$}.

\vskip.3cm
\noindent
Theorem (\ref{theorem about 3-Lie superalgebra}) can be applied to matrix Lie superalgebras. Consider a matrix Lie superalgebra $\frak{gl}(m,n)$ of $(m,n)$-supermatrices. The super trace of a supermatrix is the linear $\mathbb C$-valued functional, which satisfies the conditions of Theorem (\ref{theorem about 3-Lie superalgebra}). Indeed for any odd degree matrix $X\in\frak{gl}_1(m,n)$ it holds
$$
\mbox{Str}\,X=0,
$$
and hence the first condition for $\omega$ is satisfied. The super trace is the 1-cocycle because for any two supermatrices we have
$$
\delta\,\mbox{Str}(X,Y)=\mbox{Str}([X,Y])=0,
$$
and the second condition is also satisfied. Thus if we endow a matrix Lie superalgebra $\frak{gl}(m,n)$ with the graded triple Lie bracket
\begin{eqnarray}
[X,Y,Z]&=&\mbox{Str}(X)\,[Y,Z]+(-1)^{x\,\overline{yz}}\mbox{Str}(Y)\,[Z,X]\nonumber\\
&&\qquad\qquad\qquad\qquad +(-1)^{z\,\overline{xy}}\mbox{Str}(Z)\,[X,Y],
\label{graded triple Lie bracket for super matrices}
\end{eqnarray}
where $x$ is the degree of a super matrix $X$ and $\overline{xy}=x+y$, then $\frak{gl}(m,n)$ becomes the 3-Lie superalgebra. This matrix 3-Lie superalgebra will be denoted by $\frak{glt}(m,n)$. By analogy with the terminology of the previous section the graded triple Lie bracket (\ref{graded triple Lie bracket for super matrices}) will be referred to as \emph{the quantum super Nambu bracket}.
%%%%%%%%%%%%%%%%%%%%%%%%%%%%%%%%%%%%%%%%%%%%%%%%%%%%%%%%%%%%%%
%%%%%%%%%%%%%%%%%%%%%%%%%%%%%%%%%%%%%%%%%%%%%%%%%%%%%%%%%%%%%%
%#############################################################
%#############################################################
% NEW SECTION
\section{Basic structures of the calculus of cubic matrices}
%************************************
The aim of this section is to remind some basic notions related to the structure of a 3-dimensional matrix which we will need later in this paper to construct a 3-Lie superalgebra of cubic supermatrices. We use the terminology and the notions of the calculus of space matrices developed in \cite{Sokolov}.

\vskip.3cm
\noindent
By a 3-dimensional $m\times n\times p$-matrix $A$ we mean a system of $mnp$ complex numbers equally spaced at the points with coordinates $(x,y,z)$, where $x=0,1,\ldots,m-1, y=0,1,\ldots,n-1$ and $z=0,1,\ldots,p-1$. We will assign to each element of 3-dimensional matrix three integers $i,j,k$ and denote a corresponding element by $A_{ijk}$. If an element $A_{ijk}$ of a 3-dimensional matrix $A$ is spaced at the point with coordinates $(x,y,z)$ then integers $i,j,k$ are expressed in terms of coordinates as follow $i=x+1,\;\;j=y+1,\;\;k=z+1$. 
These relations show how we position a 3-dimensional matrix in the space relative to coordinate system. We will use the following directions in a 3-dimensional matrix: The directions $(i),(j),(k)$ of a cubic matrix coincide with the space directions of $X,Y,Z$ axes respectively. If $m=n=p$ then a 3-dimensional matrix $A$ will be called a \emph{cubic matrix} of $n$th order. The vector spaces of complex 3-dimensional $m\times n\times p$-matrices and complex cubic matrices of order $n$ will be denoted by $M^{(3)}_{m,n,p}$ and ${M}^{(3)}_n$ respectively. The vector space of rectangular complex $m\times n$-matrices will be denoted by ${M}_{m,n}$ and when $m=n$ the algebra of $n$th order square matrices will be denoted by ${M}_{n}$.

\vskip.3cm
\noindent
The elements $A_{ijk}$ of a 3-dimensional matrix $A$ with fixed first subscript $i$ are located in the plane perpendicular to $X$ axis and a collection of these elements of a 3-dimensional  matrix will be called a \emph{section of orientation} $(i)$. Similarly one can define sections of orientations $(j),(k)$. Geometrically one can get a section of a 3-dimensional matrix by cutting it by a plane perpendicular to one of the edges (or coordinate axes) of a cube. Evidently a section of any orientation can be viewed as a rectangular matrix and we will label each section of given orientation (considered as a rectangular matrix) with the help of fixed subscript. For instance, if we consider a section of orientation $(j)$ and it is determined by fixed second subscript $j=r$, where $r=1,2,\ldots,n$, then we will denote the corresponding rectangular matrix by $A^{(j)}_r$. Thus $A^{(j)}_r=(A_{irk})\in {M}_{m,p}$, where integer $r$ is fixed. A 3-dimensional matrix $A$ can be considered as the collection of $n$ rectangular $m\times p$-matrices $A^{(j)}_r$, where $r$ runs from 1 to $n$. Separating different sections by vertical lines and ordering them by section label we can map a 3-dimensional matrix $A$ to the plane of this page
\begin{equation}
\big(A^{(j)}_1\,|\,A^{(j)}_2\,|\,\ldots\,|\,A^{(j)}_n\big).
\label{collection of orientation j}
\end{equation}
Assume that $A$ is a cubic matrix of $n$th order. A section of any orientation of this matrix is a square matrix of $n$th order, and we can define the trace of a section of a cubic matrix as the trace of a square matrix. For instance given a section $A^{(j)}_r\in {M}_n$ of orientation $(j)$ we define its trace by
$$
\mbox{Tr}\,A^{(j)}_r=\sum_{i=1}^n A_{iri}.
$$

\vskip.3cm
\noindent
A collection of elements of a 3-dimensional matrix $A$ with two fixed subscripts $j,k$ is referred to as a \emph{double section of orientation} $(jk)$ or a \emph{row of direction} $(i)$. Obviously there are $np$ rows of direction $(i)$ in a 3-dimensional matrix $A$. Each row of direction $(i)$ can be considered as the row-matrix
$$
\big( A_{1jk}\;A_{2jk}\;\ldots\;A_{mjk}\big).
$$
Let $A$ be a cubic matrix of $n$th order. A collection of $n$ elements of this matrix is referred to as a \emph{transversal} of a cubic matrix if for any pair of elements of this collection there is no section of cubic matrix such that it contains these elements. For instance, if $A$ is a cubic matrix of 3rd order then three elements $A_{113}, A_{232}, A_{321}$ form the transversal of $A$. The elements of a cubic matrix equally spaced along straight lines passing through opposite vertices of a cube (diagonals of a cube) form the four transversals of a cubic matrix of $n$th order. One of these transversals $A_{111}, A_{222},\ldots, A_{nnn}$ is called the main diagonal of a cubic matrix and three others are called \emph{secondary diagonals}.

\vskip.3cm
\noindent
By analogy with the notion of transversal of a cubic matric one can define a notion of 2-dimensional transversal section of a cubic matrix. A collection of $n$ rows of direction $(i)$ is referred to as a \emph{2-dimensional transversal section of direction} $(i)$ if for any pair of rows of this collection there is no section of cubic matrix which passes through these rows. Similarly one can define a 2-dimensional transversal sections of directions $(j)$ and $(k)$. 

\vskip.3cm
\noindent
The 2-dimensional transversal section of direction $(j)$ of a cubic matrix of $n$th order
\begin{eqnarray}
(A_{111}\;\;\;\,\;A_{121}\;\,\;\ldots\;\,\;A_{1n1}),\nonumber\\
(A_{212}\;\;\;\,\;A_{222}\;\,\;\ldots\;\,\;A_{2n2}),\nonumber\\
\ldots\ldots\ldots\ldots\ldots\ldots\ldots\ldots\ldots\nonumber\\
(A_{n1n}\;\;\,A_{n2n}\;\,\ldots\;\,A_{nnn}),\nonumber
\end{eqnarray}
will be called the \emph{main diagonal section of direction} $(j)$. Analogously one can define the main diagonal sections of a cubic matrix of directions $(i)$ and $(k)$.

\vskip.3cm
\noindent
We define the \emph{Hermitian adjoint relative to direction $(j)$} matrix $A^\dag$ as the complex conjugate of a matrix $A$ and transposed with respect to the main diagonal section of direction $(j).$ Hence if we denote an element of the Hermitian adjoint matrix by $b_{ijk}$ then $b_{ijk}={\bar a}_{kji}$.

\vskip.3cm
\noindent
Similarly to the trace of a square matrix one can define a trace of a cubic matrix by means of main diagonal section \cite{Awata-Li-Minic-Yaneya}. Since there are three main diagonal sections in a cubic matrix a notion of  trace depends on a direction of a cubic matrix. Thus we define the \emph{trace relative to direction} $(j)$ of a cubic matrix $A$ as the sum of all elements of the main diagonal section of direction $(j)$. The trace relative to direction $(j)$ of a cubic matrix $A$ will be denoted by $\mbox{Tr}^{(j)}A$ and
\begin{equation}
\mbox{Tr}^{(j)}A=\sum_{i,l=1}^{n}A_{ili}.
\label{trace of cubic matrix}
\end{equation}
We can define the traces of directions $(i)$ and $(k)$ of a cubic matrix $A$ in a similar manner and we get
$$
\mbox{Tr}^{(i)}A=\sum_{k,l=1}^{n}A_{lkk},\quad
   \mbox{Tr}^{(k)}A=\sum_{j,l=1}^{n}A_{jjl}.
$$
Since a cubic matrix can be viewed as the collection of sections of chosen orientation (\ref{collection of orientation j}) it is easy to see that the trace relative to chosen direction (let say $(j)$) of a cubic matrix is equal to the sum of traces of all sections of the orientation $(j)$. For example the trace of a cubic matrix relative to the direction $(j)$ is
$$
\mbox{Tr}^{(j)}A=\sum_{r=1}^n \mbox{Tr}\,A^{(j)}_r,
$$
where $A^{(j)}_r$ is an $n$th order matrix of section of orientation $(j)$. From this it follows that $\mbox{Tr}^{(j)}A^\dag=\overline{\mbox{Tr}^{(j)}A}$.
%#########################################################################
%#########################################################################
%@@@@@@@ NEW SECTION
\section{Quantum Nambu bracket of cubic matrices}
In this section we consider a triple product of cubic matrices proposed in \cite{Awata-Li-Minic-Yaneya}. This triple product is defined with the help of the trace of a cubic matrix.

\vskip.3cm
\noindent
First we remind a notion of associativity in the case of an algebra with a ternary law of composition. Assume $\mathscr A$ is a vector space over $\mathbb C$ endowed with a triple product
$$
(a,b,c)\in \mathscr A\times\mathscr A\times\mathscr A\to (abc)\in \mathscr A.
$$
A triple product is said to be \emph{{lr}-associative} (left-right associative) \cite{Abramov2010} if it satisfies the relation
\begin{equation}
\big((abc)df\big)=\big(ab(cdf)\big).
\label{lr-associativity}
\end{equation}
A triple product is said to be an \emph{associative of the first kind} if in addition to lr-associativity (\ref{lr-associativity}) it satisfies the relation
$$
\big((abc)df\big)=\big(a(bcd)f\big).
$$
A triple product is said to be an \emph{associative of the second kind} if in addition to lr-associativity it satisfies the relation
$$
\big((abc)df\big)=\big(a(dcb)f\big).
$$

\vskip.3cm
\noindent
Let $A$ be a 3-dimensional $m\times n\times p$-matrix, $B$ be a 3-dimensional $p\times n\times q$-matrix. A \emph{product relative to direction $(j)$} of two 3-dimensional matrices $A\ast_{\small j}{\!} B$ is the 3-dimensional $m\times n\times q$-matrix whose $r$th section of orientation $(j)$ ($1\leq r\leq n$) is defined by
\begin{equation}
(A\ast_{\small j}{\!} B)^{(j)}_r=A^{(j)}_r\cdot B^{(j)}_r,
\label{binary product of 3-dim matrices}
\end{equation}
where $A^{(j)}_r$ is the $r$th section of orientation $(j)$ of a 3-dimensional matrix $A$ (considered as the $m\times p$-matrix), $B^{(j)}_r$ is the $r$th section of orientation $(j)$ of a 3-dimensional matrix $B$ (considered as the $p\times q$-matrix) and $A^{(j)}_r\cdot B^{(j)}_r$ stands for the product of two rectangular matrices. It follows from this definition that the product $A\ast B$ can be considered as the collection of $n$ sections of orientation $(j)$ (each section is an $m\times q$-matrix) and thus $A\ast_{\small j}{\!} B$ is the 3-dimensional $m\times n\times q$-matrix. If we denote this matrix by $C=A\ast_{\small j}{\!} B$ then its $r$th section of orientation $(j)$ can be expressed as follows
$$
C^{(j)}_r=A^{(j)}_r\cdot B^{(j)}_r.
$$
Similarly one can define product of 3-dimensional relative to the direction $(i)$ or $(k)$. In what follows we will use product of 3-dimensional matrices relative to direction $(j)$.

\vskip.3cm
\noindent
The rule for multiplication of two 3-dimensional matrices determined by the formula (\ref{binary product of 3-dim matrices}) can be described as follows: we take a section of orientation $(j)$ with label $r$ ($1\leq r\leq n$) in a first 3-dimensional matrix $A$ (this is a slice of a 3-dimensional matrix perpendicular to $Y$ axis), then we take a section of the same orientation $(j)$ with the same label $r$ in a second 3-dimensional matrix $B$, compute the product of these rectangular matrices and consider this product as the $r$th section of the orientation $(j)$ of the 3-dimensional matrix $C=A\ast_{j}\! B$. Virtually this multiplication of 3-dimensional matrices amounts to multiplication of rectangular matrices and consequently has properties very similar to ordinary multiplication of rectangular matrices. It is easy to see that the product defined in (\ref{binary product of 3-dim matrices}) is associative, i.e. $(A\ast_{j}\! B)\ast_{j}\! C=A\ast_{j}\! (B\ast_{j}\! C)$. Thus the multiplication of 3-dimensional matrices (\ref{binary product of 3-dim matrices}) induces the structure of unital associative algebra on the vector space of $n$th order cubic matrices $\cubic$, where the identity element of this algebra is the $n$th order cubic matrix $I$, whose any section of orientation $(j)$ is the $n$th order square unit matrix $E$, i.e. $I^{(j)}_r=E$ and $r=1,2,\ldots,n$. By other words each element of the main diagonal section of direction $(j)$ of the cubic matrix $I$ is equal to 1, $I_{iji}=1$.

\vskip.3cm
\noindent
Now $\cubic$ equipped with the binary multiplication of cubic matrices is the unital associative algebra and consequently the commutator of two cubic matrices $A,B\in \cubic$ defined as usual $[A,B]=A\ast_{\small j}{\!} B-B\ast_{\small j}{\!} A$ satisfies the Jacobi identity. Hence the unital associative algebra of cubic matrices $\cubic$ endowed with the commutator becomes the Lie algebra, which we will denote by ${\frak gl}^{(3)}_n$. It is easy to show that the trace of a cubic matrix introduced in the previous section (\ref{trace of cubic matrix}) vanishes on the commutator of two cubic matrices. Indeed we have
$$
\mbox{Tr}^{(j)}\,[A,B]=\sum_{r=1}^m \mbox{Tr}\,[A,B]^{(j)}_r
                      =\sum_{r=1}^m \mbox{Tr}\,[A^{(j)}_r,B^{(j)}_r]=0.
$$

\vskip.3cm
\noindent
In order to construct a 3-Lie algebra of cubic matrices we will need a triple product. Following \cite{Awata-Li-Minic-Yaneya} we consider the triple product of cubic matrices which is the combination of the binary product (\ref{binary product of 3-dim matrices}) of two cubic matrices and the trace of a cubic matrix
\begin{equation}
(ABC)=\mbox{Tr}^{(j)}B\,\;(A\ast_{j}\! C).
\label{triple product of cubic matrices}
\end{equation}
It is natural to call a triple product \emph{proper triple product} if it can not be reduced to a binary product, and \emph{binary generated triple product} if it is constructed by means of a binary product. Evidently the triple product (\ref{triple product of cubic matrices}) is binary generated triple product.
\begin{prop}
The triple product (\ref{triple product of cubic matrices}) of cubic matrices satisfies the following identities:
\begin{enumerate}
\item  it is lr-associative $((ABC)DF)=(AB(CDF))$,
\item $(ABC)^\dag=(C^\dag B^\dag A^\dag)$,
\item $(A(BCD)F)=(A(DCB)F)$,
\item $((ABC)DF)=((ADC)BF)$.
\end{enumerate}
\end{prop}
\noindent
We will prove only the lr-associativity of the triple product because the remained properties can be verified analogously. Making use of the definition of triple product we can write the left hand side of the lr-associativity as follows
$$
((ABC)DF)=\mbox{Tr}^{(j)}B\;\mbox{Tr}^{(j)}D\;(A\ast_{j}\! C)\ast_{j}\! F.
$$
Evidently there is the same product of scalar factors at the right hand side of the lr-associativity (because cubic matrices $B,D$ are still in the middle of the corresponding triple products), but the matrix part has the form $A\ast_{j}\! (C\ast_{j}\! F)$, which by virtue of the associativity of the binary product is equal to $(A\ast_{j}\! C)\ast_{j}\! F$.

\vskip.3cm
\noindent
Now in order to construct a 3-Lie algebra of cubic matrices we can consider the following quantum triple Nambu bracket
\begin{equation}
[A,B,C]=(ABC)+(BCA)+(CAB)-(CBA)-(BAC)-(ACB),
\label{ternary commutator}
\end{equation}
which plays the role of ternary Lie bracket for the triple product (\ref{triple product of cubic matrices}). The structure of quantum triple Nambu bracket is based on the well known representation of $S_3$, where the cyclic (even) permutations are taken with the sign "+" and the non-cyclic (odd) permutations are taken with sign "-". Evidently the quantum triple Nambu bracket (\ref{ternary commutator}) satisfies the symmetry properties of ternary Lie bracket (\ref{symmetries of ternary commutator 2}). Thus the quantum triple Nambu bracket (\ref{ternary commutator}) vanishes whenever there are two (or three) equal elements in (\ref{ternary commutator}). It is worth to mention that there is a different approach to notion of ternary commutator which is based on the group of cyclic permutations $\mathbb Z_3$ and its representation by cubic roots of unity \cite{Abramov1997}. In this approach a ternary commutator vanishes when all three elements are equal, but it is not necessarily zero if there are two equal elements.

\vskip.3cm
\noindent
Since the triple product of cubic matrices is binary generated we can write this product in terms of binary commutators as follows
\begin{equation}
[A,B,C]=\mbox{Tr}^{(j)}B\;[A,C]+\mbox{Tr}^{(j)}A\;[C,B]+\mbox{Tr}^{(j)}C\;[B,A],
\label{triple commutator expressed in binaries}
\end{equation}
and we see that the triple commutator of cubic matrices is analogous to the triple commutator of square matrices (\ref{triple commutator of square matrices}). The triple commutator (\ref{triple commutator expressed in binaries}) will be referred to as the \emph{quantum Nambu bracket for cubic matrices}. Now taking into account the property of the trace of a cubic matrix that the trace vanishes on any commutator of two cubic matrices we see that the condition (\ref{condition for omega in the theorem}) stated in Theorem (\ref{theorem 1}) is satisfied and consequently the triple commutator of cubic matrices (\ref{ternary commutator}) satisfies the Filippov-Jacobi identity. Hence the vector space of $n$th order cubic matrices $M^{(3)}_n$ equipped with the triple commutator (\ref{ternary commutator}) (or (\ref{triple commutator expressed in binaries})) is the 3-Lie algebra.

%####################################
%#################################### NEW SECTION
%************************************
\section{Cubic supermatrix and quantum super Nambu bracket}
%************************************
%************************************
%%%%%%%%%%%%%%%%%%%%%%%%%%%%%%%%%%%%%%%%%%%%%%%%%%%%%%%%%%%%%%%%%%%%%%%%%%%%%
In this section we define the $\mathbb Z_2$-graded (super) structure of a cubic matrix by assigning parities to sections of two fixed orientations and defining the parity of the intersection of two sections (which is the row of direction different from two fixed orientations) as the sum of parities of intersecting sections. Then we define the super trace of a cubic supermatrix and show that it amounts to the sum of super traces of all sections of chosen orientation.  Making use of the super trace of a supermatrix we introduce the triple product of cubic supermatrices and find the identities for this triple product. Finally we define the quantum super Nambu bracket and show that it satisfies the graded Filippov-Jacobi identity and hence induces the structure of 3-Lie superalgebra of cubic supermatrices. The residue classes modulo 2 will be denoted by $\bar 0,\bar 1$, i.e. $\mathbb Z_2=\{\bar 0,\bar 1\}$.

\vskip.3cm
\noindent
Let $A=(A_{ijk})$ be a cubic matrix of $n$th order. We define a \emph{super structure relative to the direction $(j)$} of an $n$th order cubic matrix $A$ as follows. First we fix an integer $0<r<m$, choose $r$ sections of orientation $(i)$, $r$ sections of orientation $(k)$ and assign the parity $\bar 0$ to each of chosen sections. Naturally we assign the parity $\bar 1$ to each of remained $s=n-r$ sections of orientation $(i)$ and $s$ sections of orientation $(k)$. Accordingly a section will be called even if its parity is $\bar 0$ and odd if its parity is $\bar 1$. Each row of direction $(j)$ is the intersection of one section of orientation $(i)$ and one section of orientation $(k)$. Second we define the parity of each row of direction $(j)$ as the sum of the parities of sections whose intersection is the given row of direction $(j)$. Finally we define the parity of an element of a cubic matrix to be equal to the parity of the row of direction $(j)$ which passes through this element. Similarly one can define a super structure of a cubic matrix relative to direction $(i)$ or $(k)$.

\vskip.3cm
\noindent
In the described above super structure of a cubic matrix, the choice of $r$ even sections can be made arbitrarily, but it is more convenient to order even sections, placing them in the first place.
Hence we assign the parity $\bar 0$ to first $r$ sections of orientation $(i)$, the parity $\bar 1$ to the remaining $s=n-r$ sections of the same orientation and in the same way assign the parities to the sections of orientation $(k)$. Then the sections of orientation $(i)$ with the labels $1,2,\ldots,r$ are even sections, the section of orientation $(i)$ with the labels $r+1,r+2,\ldots,m$ are odd sections and exactly the same distribution of parities takes place for sections of orientation $(k)$. A cubic matrix endowed with a super structure will be called a \emph{cubic ($r,s$)-supermatrix}. We will denote the vector space of complex $n$th order cubic ($r,s$)-supermatrices by ${M}^{(3)}_{r,s}$, where $r+s=n$.

\vskip.3cm
\noindent
Consider a section of orientation $(j)$. This section is an $n$th order square matrix denoted by $A^{(j)}_r$. Each element of this square matrix is at the intersection of given section of orientation $(j)$ and the row of direction $(j)$. But each row of direction $(j)$ has the certain parity. Consequently each element of a square matrix $A^{(j)}_r$ has the certain parity and we see that the super structure of a cubic matrix induces the super structure of each section of orientation $(j)$. Hence each $A^{(j)}_r$ is an $n$th order square supermatrix with the usual super structure of square matrices: In the upper left corner and lower right corner there are the even blocks ($r$th order square matrix and $s$th order square matrix respectively) which will be denoted by $A^{r}_{00}, A^{r}_{11}$ (here $r$ is the label of a section and we do not show the orientation $(j)$ assuming it is fixed) and two odd blocks which are the rectangular matrices $A^r_{01},A^r_{10}$ of dimensions $r\times s$ and $s\times r$ respectively. The even blocks $A^r_{00}$, where $r=1,2,\ldots,n$, form the 3-dimensional $r\times n\times r$-matrix and the even blocks $A^r_{11}$ form the 3-dimensional ($s\times n\times s$)-matrix which will be denoted by $A^{(j)}_{00},A^{(j)}_{11}$ respectively and will be called the even blocks of a cubic supermatrix $A$. Analogously we define the odd blocks of a cubic supermatrix and denote them by $A^{(j)}_{01},A^{(j)}_{10}$, where $A^{(j)}_{01}$ is the 3-dimensional $r\times n\times s$-matrix and $A^{(j)}_{10}$ is the 3-dimensional $s\times n\times r$-matrix.

\vskip.3cm
\noindent
As usual a supermatrix $A$ is said to be an \emph{even degree supermatrix} if its odd blocks $A^{(j)}_{01},A^{(j)}_{10}$ are zero 3-dimensional matrices and $A$ is said to be a \emph{odd degree supermatrix} if its even blocks $A^{(j)}_{00},A^{(j)}_{11}$ are zero 3-dimensional matrices. The degree of a supermatrix $A$ will be denoted by $|A|$ (hence $|A|$ is either $\bar 0$ or $\bar 1$) and we will also use the small letter to denote degree of a supermatrix, for example the degree of a supermatrix $A$ will be denoted by $a$. Obviously the binary product (\ref{binary product of 3-dim matrices}) of two cubic ($r,s$)-supermatrices $A,B$ has the property
$$
|A\ast_{j}\! B|=|A|+|B|,
$$
and thus
the vector space of cubic ($r,s$)-supermatrices $M^{(3)}_{r,s}$ endowed with the binary product (\ref{binary product of 3-dim matrices}) is the superalgebra.

\vskip.3cm
\noindent
Next structure on the super vector space of cubic supermatrices arises if we consider the graded binary commutator of two cubic supermatrices
\begin{equation}
[A,B]=A\ast_{j}\! B - (-1)^{ab}\;B\ast_{j}\! A,
\end{equation}
where $a,b$ are the degrees of supermatrices $A,B$ respectively. Evidently this graded commutator satisfies the graded Jacobi identity and the super vector space of cubic supermatrices becomes the Lie superalgebra, which will be denoted by ${\frak gl}^{(3)}(r,s)$.

\vskip.3cm
\noindent
Now we can extend the notion of a trace relative to direction $(j)$ of a cubic matrix to cubic supermatrices. For this purpose we will use the formula (\ref{trace of cubic matrix}). Indeed we know that each section $A^{(j)}_l$ of orientation $(j)$ is an $n$th order square ($r,s$)-supermatrix and we can compute its supertrace. Thus we define the \emph{supertrace relative to direction $(j)$ of a cubic supermatrix $A$} as follows
\begin{equation}
\mbox{Str}^{(j)} A=\sum_{l=1}^n \mbox{Str}\; A^{(j)}_l.
\end{equation}
Obviously the supertrace of any odd degree cubic supermatrix vanishes.
It is easy to verify that the supertrace defined in this way vanishes on graded commutators of cubic supermatrices. Indeed
$$
\mbox{Str}^{(j)}\;[A,B]=\sum_{l=1}^n \;\mbox{Str}\;[A^{(j)}_l,B^{(j)}_l]=0.
$$
By analogy with the triple product of cubic matrices (\ref{triple product of cubic matrices}) we define the triple product of cubic supermatrices
\begin{equation}
(ABC)=\mbox{Str}^{(j)}A\;(B\ast_j\! C).
\label{triple product of cubic supermatrices}
\end{equation}
\begin{prop}
The triple product (\ref{triple product of cubic supermatrices}) satisfies the following identities:
\begin{enumerate}
\item $(AB(CDF))=(A(CBD)F)$,
\item $((ABC)DF)=(-1)^{bc}\;((ACB)DF)$,
\item $(AB(CDF))=(CB(ADF))$,
\end{enumerate}
where $b,c$ are the degrees of cubic supermatrices $B,C$ respectively.
\end{prop}
\noindent
Next we define the graded triple commutator of cubic supermatrices by the formula
\begin{eqnarray}
[A,B,C]\!\!\!&\!=\!&\!\!\!(ABC)+(-1)^{a\,\overline{bc}}\,(BCA)+(-1)^{c\,\overline{ab}}\,(CAB)\nonumber\\
             &&-(-1)^{ab}\,(BAC)-(-1)^{bc}\,(ACB)-(-1)^{ab+bc+ac}\,(CBA),\nonumber
\end{eqnarray}
where $\overline{ab}=a+b$.
Substituting the definition of the triple product (\ref{triple product of cubic supermatrices}) into this formula we get
\begin{equation}
[A,B,C]=\mbox{Str}^{(j)}A\;[B,C]+(-1)^{a\,\overline{bc}}\mbox{Str}^{(j)}B\;[C,A]+(-1)^{c\,\overline{ab}}\mbox{Str}^{(j)}C\;[A,B].
\label{viimane valem}
\end{equation}
Since the supertrace of a cubic supermatrix vanishes on the subspace of odd degree cubic supermatrices and it also vanishes on graded commutators of cubic supermatrices we see that the both conditions of Theorem (\ref{theorem about 3-Lie superalgebra}) are satisfied and we conclude that the super vector space of cubic supermatrices equipped with the graded triple commutator (\ref{viimane valem}) is the 3-Lie superalgebra. Hence the graded triple commutator (\ref{viimane valem}) satisfies the graded Filippov-Jacobi identity and by analogy with the previously used terminology it will be referred to as the \emph{quantum super Nambu bracket for cubic supermatrices.}
%%%%%%%%%%%%%%%%%%%%%%%%%%%%%%%%%%%%%%%%%%%%%%%%%%%%%%%%%%%%%%%%%%%%%%%%%%%%%
\subsection*{Acknowledgment}
The author gratefully acknowledges that this work was financially supported by the institutional funding IUT20-57 of the Estonian Ministry of Education and Research.
%%%%%%%%%%%%%%%%%%%%%%%%%%%%%%%%%%%%%%%%%%%%%%%%%%%%%%%%%%%%%%%%%%%%%%%%%%%%%
%#############################################################################
%@@@@@@@@@@@ Bibliography

%##################################################
%##################################################
\end{document}